\documentclass{amsart}
\usepackage{latexsym,amsmath,amstext,mathrsfs,amssymb,amsthm,amscd,amsopn,verbatim,graphics,amsrefs}

\newtheorem{Thm}{Theorem}[section]
 
\newtheorem{Lem}[Thm]{Lemma}
\newtheorem{Cor}[Thm]{Corollary}

\theoremstyle{remark}
\newtheorem{Rem}[Thm]{Remark}

\theoremstyle{definition}

\newtheorem*{Reg*}{Regularization procedure}

\title[Weighted projective spaces...]{Weighted projective spaces and minimal nilpotent orbits}
\author{Carlo A. Rossi}
\address{Department of mathematics, ETH Z\"urich, 8092 Z\"urich, Switzerland}
\email{crossi@math.ethz.ch}

\begin{document}

\maketitle

\begin{abstract}
We investigate (twisted) rings of differential operators on the resolution of singularities of an irreducible component $\overline X$ of $\overline O_{\mathrm{min}}\cap\mathfrak n_+$ (where $\overline O_{\mathrm{min}}$ is the (Zarisky) closure of the minimal nilpotent orbit of $\mathfrak{sp}_{2n}$ and $\mathfrak n_+$ is the Borel subalgebra of $\mathfrak{sp}_{2n}$) using toric geometry, and show that they are homomorphic images of a certain family of associative subalgebras of $U(\mathfrak{sp}_{2n})$, which contains the maximal parabolic subalgebra $\mathfrak p$ determining $\overline O_{\mathrm{min}}$.
Further, using Fourier transforms on Weyl algebras, we show that (twisted) rings of well-suited weighted projective spaces are obtained from the same family of subalgebras.
Finally, we investigate this family of subalgebras from the representation-theoretical point of view and, among other things, re-discover in a different framework irreducible highest modules for the UEA of $\mathfrak{sp}_{2n}$. 
\end{abstract}

\section{Introduction}\label{s-I}
Let us consider the Lie algebra $\mathfrak g=\mathfrak{sp}_{2n}$: the ring $\mathcal D(\overline X)$ of differential operators on $\overline X$, an irreducible component of $\overline O_{\mathrm{min}}\cap \mathfrak n_+$, $\overline O_{\mathrm{min}}$ being the Zarisky closure of the minimal nilpotent orbit of $\mathfrak g$ and $\mathfrak n_+$ the Borel subalgebra $\mathfrak n_+$ of $\mathfrak g$, has been discussed in details in~\cite{Levasseuretal} and independently in~\cite{Musson1}; for an explicit characterization of $\overline X$, we refer to~\cite{Levasseuretal} and Section~\ref{s-1}.
The variety $\overline X$ has the property that it admits a surjective morphism
\begin{equation}\label{eq-diffop}
U(\mathfrak g)\to \mathcal D(\overline X);
\end{equation}
while Levasseur et al.~\cite{Levasseuretal} have obtained (\ref{eq-diffop}) as a special case in a more general framework, invoking Fourier transform, Musson~\cite{Musson1} highlighted the fact that $\overline X$ is a toric variety (affine and singular).
In particular, the toric structure permits~\cite{Musson2} to find an explicit realization of the epimorphism (\ref{eq-diffop}) by constructing generators and relations of $\mathcal D(\overline X)$ using homogeneous coordinates on $\overline X$.
The construction of such coordinates for $\overline X$ has the advantage of producing almost without effort an explicit resolution of singularities of $\overline X$, which we denote by $\widetilde X $: explicitly, $\widetilde X $ can be identified with the Serre bundle $\mathcal O(-2)$ over $\mathbb P^{n-1}$.
The toric structure of $\widetilde X $ allows~\cite{Musson2} to compute generators and relations for any twisted ring $\mathcal D_{\mathcal O(\ell D_0)}(\widetilde X)$, $D_0$ being a generator of the Picard group of $\widetilde X$ and $\ell$ an integer.

Our main result is a version of (\ref{eq-diffop}) for $\widetilde X $, for which we need some preparations.
Observe that the construction of $\overline X$ requires the choice of a maximal parabolic subalgebra $\mathfrak p$ of $\mathfrak g$.
In the case at hand, $\mathfrak p$ splits as $\mathfrak m\oplus \mathfrak r_+$, where $\mathfrak m$ is a reductive algebra, with one-dimensional center $\mathfrak z$ and semisimple part $\mathfrak{sl}_n$, while $\mathfrak r_+$ is the nilradical, which, as an $\mathfrak m$-module, is isomorphic to the finite-dimensional, irreducible module $L(2\varpi_1)$, $\varpi_i$, $i=1,\dots,n-1$ being the $i$-th fundamental weight of $\mathfrak{sl}_n$.
We have the triangular decomposition $\mathfrak g=\mathfrak r_-\oplus\mathfrak m\oplus\mathfrak r_+$, where $\mathfrak r_-$ is, as an $\mathfrak m$-module, isomorphic to $L(2\varpi_{n-1})$, $\varpi_i$ being the fundamental weights of $\mathfrak{sl}_n$.
\begin{Thm}
For any integer $\ell$, consider the subalgebra $\mathfrak A_\ell$ of $U(\mathfrak g)$, generated by $1$, $\mathfrak m$, $\mathfrak r_-$ and $\mathfrak r_+\mathfrak z_\ell$; then, for any integer number $\ell$, there is a surjective morphism
\[
\mathfrak A_\ell\to \mathcal D_{\mathcal O(\ell D_0)}(\widetilde X ).
\]
Here, by $\mathfrak z_\ell$ is meant the $1$-dimensional subspace of $U(\mathfrak g)$ spanned by $z_\ell:=z+\frac{\ell}2+\frac{n}4$, where $z$ is the generator of $\mathfrak z$ chosen according to Subsection 3.2, Remark $(vii)$.
\end{Thm} 
In particular, we will show that there is a surjective morphism from $U(\mathfrak p)$ to the subalgebra of $D_{\mathcal O(\ell D_0)}(\widetilde X )$ of differential operators of non-positive degree, for any integer $\ell$, and a surjective morphism from $U(\mathfrak{sl}_n)$ to the subalgebra of differential operators of degree $0$: observe that there is a natural right $SL_n$-action on $\widetilde X $, through which one obtains the vector fields inducing the differential operators in $U(\mathfrak{sl}_n)$. 
The notion of degree needs a brief explanation: all rings of differential operators considered here are subrings of the Weyl algebra $\mathcal D(\mathbb A^{n+1})$, hence all elements are linear combinations of differential monomials $Q^\mu P^\nu$, $Q_i$, resp.\ $P_i$, $i=1,\dots,n+1$, being multiplication operators by the variables $Q_i$, resp.\ differential operators w.r.t.\ $Q_i$, and $\mu$, $\nu$ are multiindices with $n+1$ components.
The degree of such a monomial is simply $\sum_{i=1}^{n+1} \tau_i$, where $\tau=\mu-\nu$.
 
As shown in~\cite{FR}, there exists an explicit isomorphism between (twisted) rings of differential operators on projective $n$-space $\mathbb P^n$ and on the blow-up $\widetilde{\mathbb A^n}$ of affine $n$-space $\mathbb A^n$: analogously, there is an isomorphism between (twisted) rings of differential operators on the weighted projective $n$-space $Y=\mathbb P^n(1,\dots,1,2)$ and on the resolution $\widetilde X $.
\begin{Thm}
For any even integer $\ell$, there is a ring isomorphism
\[
\mathcal D_{\mathcal O(\ell D_0)}(\widetilde X )\to \mathcal D_{\mathcal O(\ell-2)}(Y).
\]
\end{Thm} 
The isomorphism in the previous Theorem needs some explanations.
It is realized concretely by means of Fourier transform~\cite{MussonRueda}, using the main result of~\cite{FR}: we simply observe that the weighted projective space $Y$ and $\widetilde X $ are related to each other by means of $I$-reflections.
We notice that this phenomenon cannot be observed on $\overline X$.
Differently from usual projective space, the weighted projective space $Y$ is a singular toric variety: hence, the main Theorem of~\cite{FR} cannot be readily applied. 
Weighted projective spaces share, on the other hand, many of the properties of usual projective spaces: in particular, they possess natural Serre sheaves $\mathcal O(m)$, $m$ being an integer, which are invertible if and only if $m$ is even.
Recalling the ingredients of the main Theorem of~\cite{FR}, we have to determine for which choices of an integer $\ell$ the image of the invertible sheaf $\mathcal O(\ell D_0)$ on $\widetilde X $ in $A_{n-1}(Y)$ corresponds to an invertible sheaf on $Y$: by the above reasonings, only even integers $\ell$ will do the job.

We may further combine both Theorems to get a useful corollary.
\begin{Cor}
For any even integer $\ell$, there is a surjective morphism
\[
\mathfrak A_{\ell+2}\to \mathcal D_{\mathcal O(\ell)}(Y).
\]
\end{Cor}
As a consequence, there is a surjective morphism from $U(\mathfrak p)$ to the subalgebra $D_{\mathcal O(\ell)}(Y)$ of differential operators of non-positive degree; observe also that there is an $M$-action on $Y$, with $M$ as above.
We notice that, unlike in the case of usual projective spaces, one cannot realize weighted projective spaces as $Sp_{2n}$-homogeneous spaces.

Finally, the twisted ring $\mathcal D_{\mathcal O(\ell D_0)}(\widetilde X )$, resp.\ $\mathcal D_{\mathcal O(\ell)}(Y)$, where $\ell$ is an integer, resp.\ an even integer, acts on the cohomology groups of $\widetilde X $, resp.\ $Y$, with values in the invertible sheaf $\mathcal O(\ell D_0)$, resp.\ $\mathcal O(\ell)$.
The cohomology groups can be explicitly computed in both cases, see Subsections~\ref{ss-1-2} and~\ref{ss-3-2}, and via the previous Theorems, are endowed with the structure of modules over $\mathfrak A_\ell$ and $\mathfrak A_{\ell+2}$ respectively.
The representation theory of algebras of the form $\mathfrak A_\ell$ presents certain similar features to the representation theory of $\mathfrak g$: in particular, we may introduce a category of $\mathfrak A_\ell$-modules, which plays the r\^ole of the Bernstein--Gel'fand--Gel'fand (BGG) category $\mathcal O$ of $\mathfrak g$.
We examine this category and find a criterion expressing the possibility of such an $\mathfrak A_\ell$-module to lift to a $\mathfrak g$-module.
Of particular importance is the following result, which generalizes results of Levasseur et al~\cite{Levasseuretal} regarding the $\mathfrak g$-module of regular functions on $\overline X$.
\begin{Thm}
For any non-positive integer $\ell$, the space of global sections of the invertible sheaf $\mathcal O(\ell D_0)$ over $\widetilde X $ is isomorphic to $L(-\frac{1}2\varpi_n)$, if $\ell$ is even, or to $L(\varpi_{n-1}-\frac{3}2\varpi_n)$, if $\ell$ is odd.
\end{Thm} 
Here, $L(\mu)$, for some weight $\mu$ of $\mathfrak h$, the Cartan subalgebra of $\mathfrak g$, denotes the unique irreducible quotient of the Verma module of highest weight $\mu$; $\varpi_i$, $i=1,\dots,n$, denotes the $i$-th fundamental weight of $\mathfrak g$.
Observe that in both cases, the highest weights belong to the same Weyl orbit of weights of $\mathfrak h$.

In all other cases, we are able to characterize explicitly the cohomology groups as irreducible highest weight $\mathfrak A_\ell$- or $\mathfrak A_{\ell+2}$-modules in the aforementioned category: in particular, we are able to decompose them into irreducible, finite-dimensional $\mathfrak{sl}_n$-modules (recall the observations about the $M$-action on both varieties $\widetilde X $ and $Y$), with a grading induced by $\mathfrak z_\ell$ or $\mathfrak z_{\ell+2}$ and by the modules $\mathfrak r_+$ and $\mathfrak r_-$; see Subsubsections~\ref{sss-4-2-1} and~\ref{sss-4-2-2} for more details on cohomology groups as $\mathfrak A_\ell$-modules. 

{\bf Acknowledgements.} We thank Giovanni Felder for many useful discussions many helpful comments and Thierry Levasseur for having carefully read the first draft of this manuscript and suggested a better terminology; we are also greatly indebted with the referee of Representation Theory for many useful comments and for having pointed out a mistake and some gaps in the proof of the results of Subsection~\ref{ss-1-3}.

\section{A suitable irreducible component of $\overline O_{\mathrm{min}}\cap\mathfrak n_+$ in $\mathfrak{sp}_{2n}$}\label{s-1}
We consider the root system of $\mathfrak g=\mathfrak{sp}_{2n}$ following the notations of Bourbaki~\cite{Bourbaki}: the only simple root of $\mathfrak g$ whose coefficient in any root decomposition is at most $1$ is $\alpha_n$, hence $\alpha_n$ specifies a unique maximal parabolic subalgebra $\mathfrak p$ with decomposition into reductive part $\mathfrak m$ and nilradical $\mathfrak r_+$ containing the highest root $e_\beta$.
The irreducible component $\overline X$ of $\overline O_{\mathrm{min}}\cap\mathfrak n_+$ is, by definition, the Zarisky closure of the adjoint cone of highest weight $X=Me_\beta$, $M$ a connected, simply connected algebraic group with Lie algebra $\mathfrak m$: it corresponds to the variety of quadratic forms on $\mathbb C^n$ of rank less or equal than $1$, see e.g.\ Section 3.2 of~\cite{Levasseuretal}.
This characterization of $\overline X$ makes evident its toric structure, the corresponding fan $\Delta$ consists of a single cone in $N\cong \mathbb Z^n$ spanned by $n$ vectors $k_i$, $i=1,\dots,n$, such that $(i)$ they can be completed to $n+1$ vectors generating $N$, and $(ii)$ they satisfy $\sum_{i=1}^n k_i=2k_{n+1}$.
Observe that $\Delta$ is not regular in the sense of~\cite{FR} and $\overline X$ is singular.

For subsequent computations, we need to characterize the toric variety $\overline X$ as a quotient of a quasi-affine variety $Y$ by some torus, in the same spirit of~\cite{Musson1}.
By the previous characterization of the fan $\Delta$, the quasi-affine variety is simply the product $\mathbb A^n\times \mathbb A^1\smallsetminus\{0\}$, while the torus action is 
\[
\mathbb C^\times \times \mathbb A^n\times\mathbb A^1\smallsetminus\{0\}\ni (u,Q_1,\dots,Q_{n+1})\mapsto (uQ_1,\dots,uQ_n,u^{-2}Q_{n+1})\in \mathbb A^n\times\mathbb A^1\smallsetminus\{0\}.
\] 
\begin{Rem}\label{r-suitable}
The variety $\overline X$ is a special case ($\mathfrak g$ of type $C_n$) of many (affine and mostly singular) varieties appearing in~\cite{Levasseuretal}, for a semisimple Lie algebra $\mathfrak g$ of type $A_n$, $B_n$, $C_n$, $D_n$, $E_6$ and $E_7$: such varieties have the remarkable property that they admit a surjective homomorphism from $U(\mathfrak g)$ onto their rings of regular differential operators. 
In~\cite{Levasseuretal}, the last statement is proved in general by means of Fourier transform; in~\cite{FR}, in the particular case $\mathfrak g$ of type $A_n$, the statement is alternatively proved using partial Fourier transforms on toric varieties and classical results of~\cite{BorhoBrylinski}; finally, if $\mathfrak g$ is of type $C_n$, it is also proved in~\cite{Musson2} using toric geometry and an explicit characterization of the ring of differential operators, see also Subsection~\ref{ss-1-3} for more details. 
\end{Rem}

\subsection{Resolution of singularities}\label{ss-1-1}
We consider the fan $\widetilde \Delta$ obtained from $\Delta$, whose one-dimensional rays are spanned by the vectors $k_i$, $i=1,\dots,n+1$, and whose cones of maximal dimension are spanned by exactly $n$ generating vectors out of the $k_i$, except the cone spanned by $k_1,\dots,k_n$: $\widetilde \Delta$ is now regular in the sense of~\cite{FR} and determines a toric variety $\widetilde X $, which is, by direct construction, a resolution of singularities of $\overline X$.
Concretely, using homogeneous coordinates, $\widetilde X $ is the quotient $(\mathbb A^n\smallsetminus\{0\}\times\mathbb A^1)// \mathbb C^\times$ with torus action
\[
\mathbb C^\times \times \mathbb A^n\smallsetminus\{0\}\times\mathbb A^1\ni (u,Q_1,\dots,Q_{n+1})\mapsto (uQ_1,\dots,uQ_n,u^{-2}Q_{n+1})\in \mathbb A^n\smallsetminus\{0\}\times\mathbb A^1.
\] 
The variety $\widetilde X $ is covered by $n$ open affine subsets, each isomorphic to the affine space $\mathbb A^n$; observe that $\widetilde X $ is identified with the line bundle $\mathcal O(-2)$ over $\mathbb P^{n-1}$.
We notice that the algebraic group $GL_n$ acts from the right on $\widetilde X $ by inverting the natural $GL_n$-action in the first $n$ homogeneous coordinates of $\widetilde X $.

\subsection{The cohomology of $\widetilde X $}\label{ss-1-2}
Since $\widetilde \Delta$ is regular, the Picard group of $\widetilde X $ is isomorphic to the group $A_{n-1}(\widetilde X )$ of Weil divisors modulo linear equivalence.
More explicitly, the Picard group is spanned by the $T$-divisors $D_i$ (see Fulton~\cite{Fulton}) associated to the $1$-dimensional cones of $\widetilde \Delta$, with the additional relations $D_i\sim D_j$, $i,j\in\{1,\dots,n\},$ and $2D_i\sim -D_{n+1}$, $i\in \{1,\dots,n\}$: thus $A_{n-1}(\widetilde X )$ is spanned by, say, $D_0:=D_1$ over $\mathbb Z$.
For any Cartier divisor $\ell D_0$, where $\ell$ is an integer, the corresponding invertible sheaf $\mathcal O(\ell D_0)$ is a sheaf of covariants (see~\cite{Musson2}) associated to the character $u\mapsto u^\ell$ of $\mathbb C^\times$.
Further, the identification $\mathcal O(\ell D_0)\cong \mathcal O(\ell)$ on the divisor $D_{n+1}\cong \mathbb P^{n-1}$ holds true (one may thus view $\widetilde X $ as a ``weighted'' version of the blow-up $\widetilde{\mathbb A^n}$).
\begin{Thm}\label{t-cohres}
For any integer $\ell$, there is an isomorphism
\[
H^\bullet(\widetilde X ,\mathcal O(\ell D_0))\cong\bigoplus_{m\geq 0}H^\bullet (\mathbb P^{n-1},\mathcal O(\ell +2m)).
\]
\end{Thm}
\begin{proof}
The proof follows along the same lines of the proof of Theorem 4.7.\ in~\cite{FR}, the only difference being that the fiber coordinate on $\widetilde X $, viewed as a bundle over $\mathbb P^{n-1}$, has weight $-2$, which produces the factor $2$ in the shift on the right hand-side of the isomorphism
\end{proof}
As a consequence, the cohomology of $\widetilde X $ with values in the invertible sheaf $\mathcal O(\ell D_0)$ is non-trivial exactly in degree $0$ and $n-1$; furthermore, the $0$-th cohomology is always non-trivial and infinite-dimensional, while the $n-1$-th cohomology is finite-dimensional and non-trivial exactly when $\ell\leq -n$.

Recalling the cohomology of projective space $\mathbb P^{n-1}$, we know that the $0$-th cohomology of $\mathbb P^{n-1}$ with values in $\mathcal O(m)$ is spanned by monomials $Q_1^{\mu_1}\cdots Q_n^{\mu_n}$ with $\mu_i\geq 0$ and $\sum_{i=1}^n \mu_i=m$, $m\geq 0$, and that the $n-1$-th cohomology  is spanned (modulo coboundaries) by monomials $Q_1^{\mu_1}\cdots Q_n^{\mu_n}$ with $\mu_i<0$ and $\sum_{i=1}^n \mu_i=m$, $m\leq -n$.
Thus, concretely, in the non-trivial cases, the previous isomorphism can be written in the form
\[
H^\bullet(\mathbb P^{n-1},\mathcal O(\ell+2m))\ni Q_1^{\mu_1}\cdots Q_n^{\mu_n}\mapsto Q_1^{\mu_1}\cdots Q_n^{\mu_n}Q_{n+1}^m\in H^\bullet(\widetilde X ,\mathcal O(\ell D_0)).
\]
\subsection{Twisted rings of differential operators on $\widetilde X $}\label{ss-1-3}
Before discussing the twisted ring of differential operators $\mathcal D_{\mathcal O(\ell D_0)}(\widetilde X )$, we consider the ring $\mathcal D(\overline X)$ of global, regular differential operators on $\overline X$.
The explicit form of this ring was computed in~\cite{Musson2}; for a better understanding of forthcoming computations, we compute explicitly the generators of $\mathcal D(\overline X)$ by using the main result of~\cite{Musson1}.

First of all, by the characterization of $\overline X$ as a toric variety at the beginning of Section~\ref{s-1}, $\mathcal D(\overline X)$ is spanned by differential monomials $Q^\mu P^\nu$, where $(i)$ $\mu$ is a multiindex with $\mu_i\geq 0$, $i=1,\dots,n$, and $\mu_{n+1}\in\mathbb Z$, and $(ii)$ $\nu$ is a multiindex with non-negative components, such that $\tau:=\mu-\nu$ satisfies $\sum_{i=1}^n\tau_i-2\tau_{n+1}=0$ (we use the same notations for the Weyl algebra as in~\cite{FR}).
The degree of a differential monomial $Q^\mu P^\nu$ is defined here as $|\tau|:=\sum_{i=1}^{n+1}\tau_i$, $\tau$ as above.
\begin{Lem}\label{l-pbw-sing}
The ring $\mathcal D(\overline X)$ is spanned by $1$ and by the differential monomials
\[
Q_i P_j,\ Q_{n+1}P_{n+1},\ Q_iQ_jQ_{n+1},\ P_iP_j Q_{n+1}^{-1},\ i,j=1\dots,n,
\]
of degree $0$, $3$ and $-3$ respectively, subject to the relation $\sum_{i=1}^n Q_iP_i-2Q_{n+1}P_{n+1}=0$.
\end{Lem}
\begin{proof}
We show that any differential monomial $Q^\mu P^\nu$ as above spanning $\mathcal D(\overline X )$ can be written as a linear combination of products of the above generators; for the sake of simplicity, we adopt the notation $Q':=(Q_1,\dots,Q_n)$, and similarly for $P'$ and for multiindices $\mu'$, $\nu'$ etc.
We consider a differential monomial $Q^\mu P^\nu$ of degree $0$: then $\tau_{n+1}=0$ and $|\tau'|=0$.
Hence, we can rewrite
\[
Q^\mu P^\nu=(Q')^{\mu'}(P')^{\nu'} Q_{n+1}^{\mu_{n+1}}P_{n+1}^{\mu_{n+1}},
\]
and $(Q')^{\mu'}(P')^{\nu'}$ has degree $0$. 
Since $\tau_{n+1}=0$, the above conventions over the multiindices force $\mu_{n+1}\geq 0$.

If the degree of $Q^\mu P^\nu$ is strictly positive, then $\tau_{n+1}>0$ and $|\tau'|>0$, thus
\[
Q^\mu P^\nu=(Q')^{\lambda'}Q_{n+1}^{\tau_{n+1}} (Q')^{\mu''}(P')^{\nu''}Q_{n+1}^{\nu_{n+1}}P_{n+1}^{\nu_{n+1}}, 
\]
where $|\lambda'|=|\tau'|=2\tau_{n+1}$ and $|\tau''|=0$.
We observe that, by construction, negative powers of the component $Q_{n+1}$ do not appear here.

Finally, if $Q^\mu P^\nu$ has strictly negative degree, then $\tau_{n+1}<0$: since the component $\mu_{n+1}$ can be negative, we have to distinguish two cases, namely $(i)$ $\mu_{n+1}<0$ and $(ii)$ $\mu_{m+1}\geq 0$.
We begin by discussing $(i)$: since $\tau_{n+1}<0$, then also $|\tau'|<0$, and $Q^\mu P^\nu$ can be rewritten as 
\[
Q^\mu P^\nu=(Q')^{\mu''}(P')^{\nu''} Q_{n+1}^{\nu_{n+1}}P_{n+1}^{\nu_{n+1}} (P')^{\lambda'} Q_{n+1}^{\tau_{n+1}},
\]
where $|\lambda'|=-|\tau'|=-2\tau_{n+1}$ and $|\tau''|=0$.
As for $(ii)$, by similar arguments, and since $\mu_{n+1}\geq 0$, the differential monomial $Q^\mu P^\nu$ can be rewritten as
\[
Q^\mu P^\nu=(Q')^{\mu''}(P')^{\nu''} Q_{n+1}^{\mu_{n+1}}P_{n+1}^{\mu_{n+1}} (P')^{\lambda'} P_{n+1}^{-\tau_{n+1}},
\]
where $|\lambda'|=-|\tau'|=-2\tau_{n+1}$ and $|\tau''|=0$.

Using the commutation relations of the Weyl algebra, $Q^\mu P^\nu$ can be further rewritten in all four cases as a linear combination of products of the given generators: we only observe that, in $(ii)$, the monomial $Q^\mu P^\nu$ can be rewritten as a linear combination of products of $Q_i P_i$, $i=1,\dots,n$, $Q_{n+1}P_{n+1}$ and $P_iP_jP_{n+1}$, $i,j=1,\dots,n$. 
But the monomials $P_iP_jP_{n+1}$ equal $(P_iP_jQ_{n+1}^{-1})(Q_{n+1}P_{n+1})$, hence the claim follows.
\end{proof}
Easy computations using Lemma~\ref{l-pbw-sing} show that there is a surjective algebra homomorphism from $U(\mathfrak{sp}_{2n})$ onto $\mathcal D(\overline X)$ (see also~\cite{Levasseuretal} and~\cite{Musson2}): this surjection is explicitly given in terms of Chevalley--Cartan generators of $\mathfrak{sp}_{2n}$ by the formul\ae
\begin{equation}\label{eq-chevgen}
\begin{aligned}
e_i&=\begin{cases}
-Q_{i+1}P_i,& i=1,\dots,n-1\\
\frac{1}2 P_n^2Q_{n+1}^{-1},& i=n
\end{cases},\\ 
h_i&=\begin{cases}
-Q_iP_i+Q_{i+1}P_{i+1},& i=1,\dots,n-1\\
-Q_nP_n-\frac{1}2,& i=n
\end{cases},\\
f_i&=\begin{cases}
-Q_iP_{i+1},& i=1,\dots,n-1\\
-\frac{1}2Q_n^2Q_{n+1}
\end{cases}.
\end{aligned}
\end{equation}
We consider, for an integer $\ell$, the twisted ring of differential operators $\mathcal D^r_{\mathcal O(\ell D_0)}(\widetilde X)$ with rational coefficients: there is a well-defined homomorphism
\begin{equation}\label{eq-twist-rat}
U(\mathfrak g)\overset{\psi_\ell}\to \mathcal D^r_{\mathcal O(\ell D_0)}(\widetilde X), 
\end{equation}
which is induced by the formul\ae\ (\ref{eq-chevgen}).
\begin{Rem}\label{rem-ratop}
The ring $\mathcal D^r_{\mathcal O(\ell D_0)}(\widetilde X)$ is well-defined, as it is viewed as the ring of global sections of the twisted sheaf of differential operators on $\mathcal O(\ell D_0)$ with rational coefficients: the latter sheaf is obtained by tensoring the twisted sheaf of differential operators $\mathcal O(\ell D_0)$ by the sheaf of rational functions on $\widetilde X$ over the structure sheaf of $\widetilde X$.
The relation $\sum_{i=1}^n Q_iP_i-2 Q_{n+1}P_{n+1}-\ell=0$ in $\mathcal D_{\mathcal O(\ell D_0)}(\widetilde X)$ (see Lemma~\ref{l-pbw1}) holds true obviously also in $\mathcal D^r_{\mathcal O(\ell D_0)}(\widetilde X)$.
\end{Rem}

Similarly to Lemma~\ref{l-pbw-sing}, the toric structure of $\widetilde X$ permits to compute explicit generators for $\mathcal D_{\mathcal O(\ell D_0)}(\widetilde X)$.
\begin{Lem}\label{l-pbw1}
For any integer $\ell$, the ring $\mathcal D_{\mathcal O(\ell D_0)}(\widetilde X )$ is spanned by $1$ and by the differential monomials
\[
Q_i P_j,\ Q_{n+1}P_{n+1},\ Q_iQ_jQ_{n+1},\ P_iP_j P_{n+1},\ i,j=1\dots,n,
\]
of degree $0$, $3$ and $-3$ respectively, subject to the relation $\sum_{i=1}^n Q_iP_i-2Q_{n+1}P_{n+1}-\ell=0$.
\end{Lem}
\begin{proof}
$\mathcal D_{\mathcal O(\ell D_0)}(\widetilde X )$ is spanned by differential monomials $Q^\mu P^\nu$, where $\mu$, $\nu$ are now both multiindices with non-negative components such that $\tau=\mu-\nu$ satisfies $\sum_{i=1}^n\tau_i-2\tau_{n+1}=0$.
The proof now follows essentially along the same lines of the proof of Lemma~\ref{l-pbw-sing}, the only difference being that, when considering differential monomials of negative degree, we only have to consider the case $(ii)$, since $\mu_{n+1}\geq 0$.
\end{proof}
We recall the parabolic triangular decomposition of $\mathfrak g$, 
\[
\mathfrak g=\mathfrak p\oplus \mathfrak r_-=\mathfrak r_+\oplus\mathfrak m\oplus \mathfrak r_-,
\]
where $\mathfrak m\cong \mathfrak{gl}_n$ is the reductive part of $\mathfrak p$, with semisimple part $\mathfrak{sl}_n$ and $1$-dimensional center $\mathfrak z$; $\mathfrak r_+$ is the abelian nilradical, and $\mathfrak r_-$ is abelian and an $\mathfrak m$-module.
In terms of differential operators, direct computations imply that the homomorphic images of $\mathfrak m$, resp.\ $\mathfrak r_-$, w.r.t.\ (\ref{eq-twist-rat}) are spanned by the monomials $Q_i P_j$, $i,j=1,\dots,n$ and $Q_{n+1}P_{n+1}$, resp.\ $Q_iQ_jQ_{n+1}$, $i,j=1,\dots,n$, with the additional relation $\sum_{i=1}^n Q_iP_i-2 Q_{n+1}P_{n+1}-\ell=0$: obviously, these two subspaces lie in $\mathcal D_{\mathcal O(\ell D_0)}(\widetilde X)$.
On the other hand, e.g.\ by Subsection 3.2, Remark $(vii)$ of~\cite{Levasseuretal}, there is a canonical generator $z$ of $\mathfrak z$, which satisfies the commutation relations
\[
[z,x]=x,\quad x\in\mathfrak r_+,\quad [z,y]=-y,\quad y\in\mathfrak r_-.
\]
The image of the generator $z$ w.r.t.\ (\ref{eq-twist-rat}) is 
\[
\psi_\ell(z)=-\frac{1}2\sum_{i=1}^n Q_iP_i-\frac{n}4=-Q_{n+1}P_{n+1}-\frac{\ell}2-\frac{n}4
\]
in $\mathcal D_{\mathcal O(\ell D_0)}(\widetilde X)$.
We denote from now on by $z_\ell$ the central element $z+\frac{\ell}2+\frac{n}4$ of $U(\mathfrak g)$.
We finally consider the differential monomials $P_iP_jP_{n+1}$ in $\mathcal D^r_{\mathcal O(\ell D_0)}(\widetilde X)$: the relation
\[
P_iP_jP_{n+1}=(-P_iP_jQ_{n+1}^{-1})(-Q_{n+1}P_{n+1})
\]
holds true.
Direct computations using Lemma~\ref{l-pbw-sing} show that the differential monomials $P_iP_jQ_{n+1}^{-1}$ lie in the homomorphic image of $\mathfrak r_+$ w.r.t.\ (\ref{eq-twist-rat}), while $Q_{n+1}P_{n+1}$ lies in the homomorphic image of $\mathfrak z\oplus\mathbb C$: it follows that the homomorphic image w.r.t.\ (\ref{eq-twist-rat}) of the subspace $\mathfrak r_+^\ell$ of $U(\mathfrak g)$ spanned by elements of the form $x z_\ell$, where $x$ is in $\mathfrak r_+$, lies in $\mathcal D_{\mathcal O(\ell D_0)}(\widetilde X)$ by Lemma~\ref{l-pbw1} and corresponds to the vector subspace spanned by the differential monomials $P_iP_jP_{n+1}$.

Thus, combining these arguments with Lemma~\ref{l-pbw1}, the homomorphism (\ref{eq-twist-rat}), restricted to the associative subalgebra $\mathfrak A_\ell$ of $U(\mathfrak g)$ generated by $1$, $\mathfrak m$, $\mathfrak r_-$ and $\mathfrak r_+^\ell$ descends to a surjective homomorphism onto $\mathcal D_{\mathcal O(\ell D_0)}(\widetilde X)$.
\begin{Thm}\label{t-parab}
For any integer $\ell$, there exists an associative subalgebra $\mathfrak A_\ell$ of $U(\mathfrak g)$ and a surjective homomorphism from $\mathfrak A_\ell$ to the twisted algebra $\mathcal D_{\mathcal O(\ell D_0)}(\widetilde X )$.
\end{Thm}

\section{Representation theory of the algebra $\mathfrak A_\ell$}\label{s-2}
For a given integer $\ell$, the algebra $\mathfrak A_\ell$, defined at the end of Subsection~\ref{ss-1-3}, contains a copy of the semisimple part of $\mathfrak m$, which is simply $\mathfrak{sl}_n$.
We have thus the additional relation in $\mathfrak A_\ell$
\[
[z_\ell,x z_\ell]=xz_\ell,
\]
for any $x\in\mathfrak r_+$. 
We observe that the nilradical $\mathfrak r_+$ is isomorphic to the finite-dimensional, irreducible highest weight $\mathfrak{sl}_n$-module $L(2\varpi_1)$, with highest weight vector $e_\beta$ and highest weight $\beta$, $\beta$ being the highest root of $\mathfrak g$; moreover, $\mathfrak r_+\cong\mathfrak r_-$ as a vector space, while, as an $\mathfrak m$-module, it is isomorphic to the finite-dimensional, irreducible highest weight module $L(2\varpi_{n-1})$.   
On the other hand, using again the Cartan involution on $\mathfrak{sl}_n$, $\mathfrak r_+$ can be viewed as the finite-dimensional, irreducible lowest weight module $L(-2\varpi_{n-1})$, with lowest weight vector $e_{\alpha_n}$ and lowest weight $\alpha_n$; similar results hold true for $\mathfrak r_-$.

\subsection{A natural ``universal'' parabolic subalgebra of $\mathfrak A_\ell$}\label{ss-2-1}
We consider the vector subspace $\mathfrak p_\ell$ of $\mathfrak A_\ell\subset U(\mathfrak g)$ spanned by $\mathfrak m$ and $\mathfrak r_+^\ell$.
\begin{Lem}\label{l-parab}
The vector space $\mathfrak p_\ell$ inherits from $U(\mathfrak g)$ the structure of a Lie algebra; furthermore, as a Lie algebra, $\mathfrak p_\ell$ is isomorphic to the parabolic subalgebra $\mathfrak p$ of $\mathfrak g$.
\end{Lem}
\begin{proof}
The Lie algebra structure is inherited from the associative structure of $U(\mathfrak g)$: since $\mathfrak r_+$ is an $\mathfrak m$-module and $z_\ell$ commutes in $U(\mathfrak g)$ with $\mathfrak m$, the non-trivial statement to prove is that $\mathfrak r_+ z_\ell$ is abelian.
We consider $x$, $x'$ in $\mathfrak r_+$: hence
\begin{align*}
[x z_\ell,x' z_\ell]&=[x z_\ell,x']z_\ell+x'[x z_\ell,z_\ell]=\\
&=x[z_\ell,x']z_\ell+[x,x']z_\ell^2+x'x[z_\ell,z_\ell]+x'[x,z_\ell]z_\ell=\\
&=xx'z_\ell-x'xz_\ell=[x,x']z_\ell=0,
\end{align*}
where we used that $\mathfrak r_+$ is an abelian Lie algebra.
The Lie algebra isomorphism $\mathfrak p\cong \mathfrak p_\ell$ is simply given by 
\[
\mathfrak m\ni m\mapsto m\in\mathfrak p_\ell,\quad r_+\ni x\mapsto x z_\ell.
\]
\end{proof}
As a result, for any integer $\ell$, the algebra $\mathfrak A_\ell$ contains a copy of the parabolic subalgebra $\mathfrak p\subset \mathfrak g$; filtration arguments \`a la Poincar\'e--Birkhoff--Witt show that this induces an isomorphism between $U(\mathfrak p)$ and the subalgebra of $\mathfrak A_\ell$ spanned by $\mathfrak m$ and $\mathfrak r_+^\ell$.
We observe that, in spite of Theorem~\ref{t-parab} of Subsection~\ref{ss-1-3}, $\mathfrak p\cong\mathfrak p_\ell$ is the vector space spanned by $Q_i P_j$, $Q_{n+1}P_{n+1}$ and $P_iP_jP_{n+1}$, $i,j=1,\dots,n$.

\subsection{The category $\mathcal O_{\mathfrak A_\ell}$}\label{ss-2-2}
We consider an $\mathfrak A_\ell$-module $M$: by previous considerations, $M$ inherits the structure of a $U(\mathfrak{sl}_n)$-module.
An element of the dual of $\mathfrak z_\ell$, the $1$-dimensional subspace spanned by $z_\ell$, is said to be a weight for $\mathfrak A_\ell$: it is simply a complex number.
$M$ is said to be $\mathfrak z_\ell$-diagonalizable, if it splits into a direct sum 
\[
M=\bigoplus_{\lambda\in\mathfrak z_\ell^*} M^\lambda,\quad M^\lambda=\left\{v\in M: z_\ell v=\lambda v\right\}.
\]
A weight of $M$ is a weight $\lambda$ with non-trivial weight subspace $M^\lambda$.
We introduce a partial order on $\mathfrak z^*_\ell$ via 
\[
\mu\leq \nu\Leftrightarrow \nu-\mu\in \mathbb N.
\] 
Previous computations yield 
\[
(xz_\ell)(M^\lambda)\subseteq M^{\lambda+1},\quad y(M^\lambda)\subseteq M^{\lambda-1},\quad m(M^\lambda)\subseteq M^\lambda,
\]
for $x\in\mathfrak r_+$, $y\in\mathfrak r_-$ and $m\in\mathfrak m$.

We define the category $\mathcal O_{\mathfrak A_\ell}$ by the following requirements: an object of $\mathcal O_{\mathfrak A_\ell}$ is an $\mathfrak A_\ell$-module such that
\begin{enumerate}
\item[$i)$] $M$ is a $\mathfrak z_\ell$-diagonalizable $\mathfrak A_\ell$-module;
\item[$ii)$] every weight subspace $M^\lambda$ of $M$ is finite-dimensional;
\item[$iii)$] there exists a weight $\mu$, such that all weights $\lambda$ of $M$ satisfy $\lambda\leq \mu$; hence, $\mu$ is called a highest weight.
\end{enumerate}
By standard arguments about gradations, the category $\mathcal O_{\mathfrak A_\ell}$ is closed w.r.t.\ taking submodules and quotient modules.

Observe that for the weight subspace to the highest weight $\mu$ of an object $M$ of $\mathcal O_{\mathfrak A_\ell}$, we have $(\mathfrak r_+^\ell)M^\mu=0$. 
Furthermore, every weight subspace $M^\lambda$ of $M$ in $\mathcal O_{\mathfrak A_\ell}$ has the structure of an $\mathfrak{sl}_n$-module.

Since every weight subspace $M^\lambda$ of $M$ is finite-dimensional, it splits into a finite direct sum of irreducible, finite-dimensional $\mathfrak{sl}_n$-modules $L(\nu):=L_{\mathfrak{sl}_n}(\nu)$, for $\nu$ a dominant weight, possibly depending on the $\mathfrak z$-weight $\lambda$.
In particular, every weight subspace $M^\lambda$ of an object $M$ of $\mathcal O_{\mathfrak A_\ell}$ contains finitely many highest weight vectors w.r.t.\ the $\mathfrak{sl}_n$-action; if, moreover, $\lambda$ is a highest weight w.r.t.\ $\mathfrak z_\ell$, such vectors are also annihilated by $\mathfrak r_+^\ell$.
As a consequence, an object $M$ of $\mathcal O_{\mathfrak A_\ell}$ possesses finitely many highest weight vectors $v$ for $\mathfrak A_\ell$, i.e.\ vectors $v$ in $M$ satisfying
\[
e_{\alpha_i}v=0,\ i=1,\dots,n-1,\ (e_{\alpha_n}z_\ell)v=0,\ z_\ell v=\mu z,\ hv=\nu(h)v,\ 
\]
where $h$ belongs to the Cartan subalgebra $\mathfrak h$ of $\mathfrak{sl}_n$, and $\nu$ is a weight for $\mathfrak h$.
\begin{Lem}\label{l-hw}
If $v$ is a highest weight vector of $M$ in $\mathcal O_{\mathfrak A_\ell}$, the vectors $e_{\alpha_n}^m v$, for any positive integer $m$, constitute a sequence of $\mathfrak{sl}_n$-primitive weight vectors in $M^{\mu-m}$, of dominant weight $\nu-m\alpha_n|_{\mathfrak h}$. 
\end{Lem}
\begin{proof}
The proof follows by a standard induction argument; the main point is that any bracket $[e_{\alpha_i},e_{-\alpha_n}]$, $i=1,\dots,n-1$, vanishes, since $-\alpha_n+\alpha_i$ does not belong to the root system of $\mathfrak g$.
Notice that $\alpha_n|_\mathfrak h=-2\varpi_{n-1}$, hence, if $\nu$ is a dominant weight for $\mathfrak{sl}_n$, so is also $\nu-m\alpha_n|_\mathfrak h$, for any positive integer $m$.
\end{proof}
(We observe that the vectors $e_{-\alpha_n}^mv$ can also be trivial.)

Special objects of the category $\mathcal O_{\mathfrak A_\ell}$ are modules, for which every weight subspace $M^\lambda$ is irreducible as an $\mathfrak{sl}_n$-module: as a consequence, to every weight subspace $M^\lambda$ belongs a unique dominant weight $\nu$ (which possibly depends on $\lambda$), such that $M^\lambda=L(\nu)$.
If e.g.\ $v$ is a highest weight vector of such an object $M$ of $\mathcal O_{\mathfrak A_\ell}$, and if $e_{-\alpha_n}^mv\neq 0$, for any positive integer $m$, then we have
\[
M=\bigoplus_{\lambda\leq \mu}M^\lambda,\quad M^\lambda\cong L(\nu-(\mu-\lambda)\alpha_n|_\mathfrak h),
\]
which will happen for most of the examples discussed in Subsection~\ref{ss-4-2}.
Since the category $\mathcal O_{\mathfrak A_\ell}$ is closed w.r.t.\ taking submodules, a non-trivial submodule $N$ of $M$ in $\mathcal O_{\mathfrak A_\ell}$ possesses a primitive subspace, i.e.\ there is a weight $\lambda$ of $M$, such that $N^\lambda=N\cap M^\lambda$ is annihilated by $\mathfrak r_+^\ell$; obviously, $N^\lambda$ is an $\mathfrak{sl}_n$-submodule of $M^\lambda$. 
On the other hand, if for an object $M$ of $\mathcal O_{\mathfrak A_\ell}$, there exists a weight $\lambda$ and an $\mathfrak{sl}_n$-submodule $N^\lambda\subseteq M^\lambda$, which is annihilated by $\mathfrak r_+^\ell$, then $M$ is reducible: we consider the $\mathfrak A_\ell$-module $N$ generated by $N^\lambda$, which is obviously a submodule of $M$ and belongs to the category $\mathcal O_{\mathfrak A_\ell}$.
Therefore, the irreducibility of objects of $\mathcal O_{\mathfrak A_\ell}$ is in one-to-one correspondence with the existence of primitive weight subspaces, i.e.\ $\mathfrak{sl}_n$-submodules of some weight subspaces, which are annihilated by the action of $\mathfrak r_+^\ell$. 

Next, the category $\mathcal O_{\mathfrak A_\ell}$ is ``too big'': in fact, as the next Theorem shows, many objects of $\mathcal O_{\mathfrak A_\ell}$ can be regarded as $\mathfrak g$-modules.
\begin{Thm}\label{t-lift}
If the highest weight $\mu$ in the weight decomposition of an object $M$ of the category $\mathcal O_{\mathfrak A_\ell}$ does not belong to $\mathbb N$, then $M$ lifts to a $\mathfrak g$-module.
\end{Thm}
\begin{proof}
Every weight $\lambda$ in the weight decomposition of $M$ is non-zero, if $\mu$ does not belong to $\mathbb N$, since the weights of $M$ are all less or equal than $\mu$ w.r.t.\ the above partial order.
We define a $\mathfrak g$-action on $M$ by defining it on any weight subspace $M^\lambda$ of $M$:
\[
xv:=\frac{1}\lambda (xz_\ell)v,\ mv:=mv,\ yv:=yv,\ x\in\mathfrak r_+,\ m\in\mathfrak m,\ y\in\mathfrak r_-,\ v\in M^\lambda.
\]
We have to show that the previous formul\ae\ define a true action: the only non-trivial relations to prove are
\[
[m,x]v=m(xv)-x(mv)\quad \text{and}\quad [x,y]v=x(yv)-y(xv),
\]
for any $x\in\mathfrak r_+$, $m\in\mathfrak m$, $y\in\mathfrak r_-$ and $v\in M^\lambda$, and any weight $\lambda$ in the weight decomposition.
Since $\mathfrak r_+^\ell$ is an $\mathfrak m$-module by Lemma\ref{l-parab}, then 
\begin{align*}
[m,x]v&=\frac{1}\lambda ([m,x]z_\ell)v=\frac{1}\lambda [m,xz_\ell]v=\\
&=\frac{1}\lambda m((xz_\ell)v)-\frac{1}\lambda (xz_\ell)(mv)=\\
&=m(xv)-x(mv),
\end{align*}
since the action of $\mathfrak m$ preserves the weight $\lambda$, and by the obvious relation $[m,xz_\ell]=[m,x]z_\ell$ in $\mathfrak A_\ell$.
As for the second relation, we have
\begin{align*}
x(yv)-y(xv)&=\frac{1}{\lambda-1} (xz_\ell)(yv)-\frac{1}\lambda y((xz_\ell)v)=\\
&=\frac{1}{\lambda(\lambda-1)}\left(\lambda (xz_\ell)(yv)-(\lambda-1)y((xz_\ell)v)\right)=\\
&=\frac{1}{\lambda(\lambda-1)}\left(\lambda (xz_\ell)(yv)-(\lambda-1)[y,xz_\ell]v-(\lambda-1)(xz_\ell)(yv)\right)=\\
&=\frac{1}{\lambda(\lambda-1)}\left(-(\lambda-1)[y,xz_\ell]v+(xz_\ell)(yv)\right).
\end{align*}
Since $\mathfrak A_\ell$ is a subalgebra of $U(\mathfrak g)$, we have the relations in $\mathfrak A_\ell$:
\[
[y,xz_\ell]=-[x,y]z_\ell+xy,\ (z_\ell-1)(xy)=(xz_\ell)y;
\]
the latter can be also rewritten the form 
\[
(xy)v=\frac{1}{\lambda-1}(xz_\ell)(yv),\quad v\in M^\lambda,\quad \lambda,\lambda-1\neq 0.
\]
Finally, using these relations, we have
\begin{align*}
&\phantom{=}\frac{1}{\lambda(\lambda-1)}\left(-(\lambda-1)[y,xz_\ell]v+(xz_\ell)(yv)\right)=\\
&=\frac{1}\lambda [x,y](z_\ell v)-\frac{1}\lambda (xy)v+\frac{1}{\lambda(\lambda-1)}(xz_\ell)(yv)=[x,y]v.
\end{align*}
\end{proof}
Thus, $\mathfrak A_\ell$-modules, which can truly belong to $\mathcal O_{\mathfrak A_\ell}$, are modules, whose highest weight $\mu$ w.r.t.\ $\mathfrak z_\ell$ is a positive integer.
If $M$ is such a module, the sequence of weights of $M$ is contained in an arithmetic sequence of the form $\mu-n$, $n\in\mathbb N$, with $\mu$ a positive integer; as we will see in forthcoming examples, the sequence of weights can be infinite or can contain only finitely many terms.

Assume finally that every weight subspace $M^\lambda$ of such a module $M$ is an irreducible $\mathfrak{sl}_n$-module and that the highest weight $\mu$ of $M$ is non-negative: then the maximal non-trivial submodule of $M$ is associated to the maximal weight $\lambda$ of $M$, such that $M^\lambda$ is annihilated by $\mathfrak r_+^\ell$.
Namely, any submodule $N$ of $M$ is associated to a non-trivial $\mathfrak{sl}_n$-submodule of $M^\lambda$, for some weight $\lambda$ of $M$, which coincides with $M^\lambda$, since $M^\lambda$ is irreducible: thus, the corresponding submodule $N$ is $\bigoplus_{\lambda'\leq\lambda} M^{\lambda'}$.
This fact implies that, choosing $\lambda$ to be maximal among all weights of $M$ such that $M^\lambda$ is annihilated by $\mathfrak r_+^\ell$, the submodule $N:=\bigoplus_{\lambda'\leq\lambda} M^{\lambda'}$ is the maximal non-trivial submodule of $M$.
Therefore, the quotient module $M/N$ is the unique irreducible finite-dimensional quotient module of $M$ of highest weight $\mu\in\mathbb N$ w.r.t.\ $\mathfrak z_\ell$: to the highest weight $\mu$ corresponds a unique dominant highest weight $\nu$ for $M^\mu$ as a $\mathfrak{sl}_n$-module, such that $M^\mu=L(\nu)$.

\section{Weighted projective spaces}\label{s-3}
\subsection{Weighted projective spaces via $I$-reflections}\label{ss-3-1}
We consider the regular fan $\widetilde \Delta$ of the resolution of singularities $\widetilde X $ and the $I$-reflection (see~\cite{FR} and~\cite{MussonRueda})
\[
k_i\mapsto \begin{cases}
k_i, & i=1,\dots,n\\
-k_i,& i=n+1.
\end{cases}
\] 
This $I$-reflection determines generating vectors $k_i'$ and a corresponding fan $\Delta'$, whose cones of maximal dimension are spanned by exactly $n$ out of the $k_i'$.
Observe that it is not anymore a regular fan, since the cone associated to $k_1',\dots,k_n'$ is not spanned by (a part of) a basis of $N$; the generating vectors $k_i'$ satisfy the relation $\sum_{i=1}^n k_i'+2k_{n+1}'=0$.
Following~\cite{Musson2}, the toric variety $Y$ associated to $\Delta'$ is the algebro-geometric quotient $\mathbb A^{n+1}\smallsetminus \{0\}/\mathbb C^\times$, with torus action 
\[
\mathbb C^\times \times \mathbb A^{n+1}\smallsetminus \{0\}\ni (u,Q_1,\dots,Q_{n+1})\mapsto (uQ_1,\dots, uQ_n,u^2Q_{n+1})\in \mathbb A^{n+1}\smallsetminus \{0\}.
\]
Therefore, the $n$-dimensional weighted projective space $Y=\mathbb P^n(1,\dots,1,2)$ is related to the resolution of singularities $\widetilde X $ of the irreducible component $\overline X$ of $\overline O_{\mathrm{min}}\cap\mathfrak n_+$ by an $I$-reflection.

\subsection{Cohomology of weighted projective spaces}\label{ss-3-2}
The weighted projective space $Y$ can be viewed as the projective spectrum of the weighted polynomial ring $S=\mathbb C[Q_1,\dots,Q_{n+1}]$, with degrees $\mathrm{deg}\ Q_i=1$, if $i=1,\dots,n$, and $\mathrm{deg}\ Q_{n+1}=2$; thus, the component $S_d$ of degree $d$ is spanned by monomials $Q^\mu$, with the multiindex $\mu$ satisfying $\sum_{i=1}^n \mu_i+2\mu_{n+1}=d$.
By construction, $Y$ is covered by $n+1$ open affine subsets, corresponding to $\{Q_i\neq 0\}$: the affine subsets corresponding to $\{Q_i\neq 0\}$, $i=1,\dots,n$, are isomorphic to the affine space $\mathbb A^n$, and the subset $\{Q_{n+1}\neq 0\}$ is isomorphic to $\overline X$ as in Subsection~\ref{ss-1-1}.

In complete analogy with the theory of usual projective spaces, the Serre sheaf $\mathcal O(\ell)$ on $Y$, for any integer $\ell$, is the sheaf associated to the shifted graded ring $S[\ell]$ (\cite{Hartsh}, Proposition 5.11 and Definition right before).
The main difference with the usual projective space is that $\mathcal O(\ell)$ is invertible exactly when $\ell$ is even, since $S$ is generated over $S_0=\mathbb C$ by $S_1$ and $S_2$.
Using the toric structure of $Y$, we can compute the Picard group of $Y$ and $A_{n-1}(Y)$: the group $A_{n-1}(Y)$ of Weil divisors on $X$ modulo linear equivalence is isomorphic to $\mathbb Z$ and is generated, say, by the divisor $D_1$. 
More precisely, $A_{n-1}(Y)$ is spanned over $\mathbb Z$ by the Weil divisors $D_i$, subject to the relations $D_i\sim D_j$, $i,j=1,\dots,n$, and $2D_i\sim D_{n+1}$, $i=1,\dots,n$. 
On the other hand, by direct computations involving Cartier divisors, the Picard group of $Y$ is generated over $\mathbb{Z}$ by the Weil divisor $2 D_1$.
Observe the obvious identification $\mathcal O(\ell D_0)=\mathcal O(\ell)$, for any even integer $\ell$: hence, the Serre bundle $\mathcal O(2)$ is the generator of the Picard group of $Y$. 
\begin{Thm}\label{thm-coh-proj}
The cohomology of $Y$ with values in any sheaf $\mathcal O(\ell)$, for any integer $\ell$, is concentrated in degree $0$ and $n$; more explicitly, we have
\[
\bigoplus_{\ell\in\mathbb Z}H^0(Y,\mathcal O(\ell))\cong S,\quad H^n(Y,\mathcal O(-n-2))\cong\mathbb C, 
\]
and there is a perfect pairing 
\[
H^n(Y,\mathcal O(\ell))\otimes H^0(Y,\mathcal O(-\ell-n-2))\to H^n(Y,\mathcal O(-n-2))\cong \mathbb C.
\]
\end{Thm}   
\begin{proof}
The proof goes along the same lines as the proof of Theorem 5.1, pages 225--228, of~\cite{Hartsh}: the only difference is that we have to keep track of different degrees in the polynomial ring defining $Y$, which produces different shifts for $\ell$ in the above pairing.
\end{proof}
As a consequence, $H^0(X,\mathcal O(\ell))$ is non-trivial exactly for $\ell\geq 0$; on the other hand, $H^n(X,\mathcal O(\ell))$ is spanned by negative monomials $Q_1^{\mu_1}\cdots Q_{n+1}^{\mu_{n+1}}$, where $\mu_i<0$ and $\sum_{i=1}^n\mu_i+2\mu_{n+1}=\ell$ (modulo coboundaries), for $\ell\leq -n-2$.
The generator of $H^n(X,\mathcal O(-n-2))$ is thus $Q_1^{-1}\cdots Q_{n+1}^{-1}$, and the perfect pairing is given by the same formula as in the proof of Theorem 5.1 of~\cite{Hartsh}.

\section{Fourier transforms and highest weight modules}\label{s-4}
\subsection{Fourier transforms, weighted projective spaces and resolution of singularities}\label{ss-4-1}
We refer to Section 3 of~\cite{FR} for notations and for the main Theorem.
From Section~\ref{s-3}, the fan of the weighted projective space $Y$ is not regular: hence, Theorem 3.2 of~\cite{FR} requires more care, as observed in Remark 3.3 of~\cite{FR}.
Namely, we have to determine the Cartier divisors on $\widetilde X $, whose image w.r.t.\ $\phi_I$ ($\phi_I$ is defined in~\cite{FR}, Lemma 3.1) are Cartier divisors in $Y$.
As seen before, the Picard group of $\widetilde X $ is generated over $\mathbb Z$ by the divisor $D_0$, and the image w.r.t.\ $\phi_I$ of a general Cartier divisor $\ell D_0$, for an integer $\ell$, is readily computed
\[
\phi_I(\ell D_0)=\ell D_1'-D_{n+1}'\sim (\ell-2) D_1',
\]  
where $D_i$ and $D_i'$ are the $T$-Weil divisors of $\widetilde X $ and $Y$ respectively; by results of Subsection~\ref{ss-3-2}, $\ell$ must be even.
Thus, Theorem 3.2 of~\cite{FR} yields the following
\begin{Thm}\label{t-weightres}
Using the notations of Section 3 of~\cite{FR}, for any even integer $\ell$, the Fourier transform $F_I$ determines an isomorphism
\[
\mathcal D_{\mathcal O(\ell D_0)}(\widetilde X ) \to\mathcal D_{\mathcal O(\ell-2)}(Y).
\]
\end{Thm}  
\begin{Rem}\label{rem-VdB}
A partial version of the isomorphism of Theorem~\ref{t-weightres} (without geometric insights) appears in~\cite{VdB}, Corollary 2.3, where the author discusses similar rings of differential operators.
\end{Rem}
Combining Theorem~\ref{t-weightres} with the results of Subsection~\ref{ss-1-2} gives
\begin{Cor}\label{cor-sympl}
For any even integer $\ell$, there is a surjective algebra homomorphism 
\[
\mathfrak A_{\ell+2}\to \mathcal D_{\mathcal O(\ell)}(Y).
\]
\end{Cor}

\subsection{The module structure on cohomology}\label{ss-4-2}
By Theorem~\ref{t-parab} from Subsection~\ref{ss-1-2}, the cohomology groups of $\widetilde X $ with values in $\mathcal O(\ell D_0)$, for any integer $\ell$, are $\mathfrak A_\ell$-modules.
Hence, Corollary~\ref{cor-sympl} endows the cohomology groups of $Y$ with values in $\mathcal O(\ell)$, for an even integer $\ell$, with the structure of $\mathfrak A_{\ell+2}$-modules.
We now proceed to analyze the respective module structures with the tools of Section~\ref{s-2}.
\begin{Rem}\label{rem-notations}
In the forthcoming Theorems~\ref{t-0res-pos},~\ref{t-0res-neg} and~\ref{t-hres} of Subsubsection~\ref{sss-4-2-1}, and Theorems~\ref{t-0weight} and~\ref{t-nweight} of Subsubsection~\ref{sss-4-2-2}, the notation $L(\bullet)$ for the unique irreducible subquotient of a certain Verma module appears frequently, and in certain cases it refers to $\mathfrak{sl}_n$, while in some other cases to $\mathfrak{sp}_{2n}$: we will use a subscript with the corresponding Lie algebra in order to avoid confusion.
\end{Rem}

\subsubsection{The module structure on the cohomology of $\widetilde X $}\label{sss-4-2-1}
The $0$-th cohomology group $H^0(\widetilde X ,\mathcal O(\ell D_0))$ is always non-trivial and infinite-dimensional; it is spanned by monomials $Q^\mu$, for multiindices $\mu$ satisfying $\sum_{i=1}^n \mu_i-2\mu_{n+1}=\ell$.
The $n-1$-th cohomology group $H^{n-1}(\widetilde X ,\mathcal O(\ell D_0))$ is non-trivial exactly for $\ell\leq -n$: it is spanned (modulo coboundaries) by monomials $Q^\mu$, where $\mu_i<0$, $i=1,\dots,n$, $\mu_{n+1}\geq 0$, and $\sum_{i=1}^n \mu_i-2\mu_{n+1}=\ell$.
\begin{Thm}\label{t-0res-pos}
For any positive integer $\ell$, there is a weight space decomposition
\[
H^0(\widetilde X ,\mathcal O(\ell D_0))\cong \bigoplus_{\lambda\leq 0} L_{\mathfrak{sl}_n}((\ell-2\lambda)\varpi_{n-1})
\]
in the category $\mathcal O_{\mathfrak A_\ell}$; the module $H^0(\widetilde X ,\mathcal O(\ell D_0))$ is irreducible and is generated by the highest weight vector of the subspace of weight $0$ w.r.t.\ $\mathfrak z_\ell$.
\end{Thm}
\begin{proof}
By results of Subsection~\ref{ss-2-2}, we need to show that $M=H^0(\widetilde X ,\mathcal O(\ell D_0))$ is in the category $\mathcal O_{\mathfrak A_\ell}$.
The generator of $\mathfrak z_\ell$ corresponds to $-Q_{n+1}P_{n+1}$: it is obvious that the decomposition from Theorem~\ref{t-cohres} corresponds to the weight decomposition 
\[
H^0(\widetilde X ,\mathcal O(\ell D_0))=\bigoplus_{\lambda\leq 0}M^\lambda,\quad M^\lambda=H^0(\mathbb P^{n-1},\mathcal O(\ell-2\lambda)).
\]
The weight subspaces $M^\lambda$ are finite-dimensional; the sequence of weights coincides with the non-positive integers.
Observe that the highest weight w.r.t.\ $\mathfrak z_\ell$ belongs to $\mathbb N$.
Every weight subspace $M^\lambda$ is a finite-dimensional irreducible $\mathfrak{sl}_n$-module, with highest weight vector $Q_n^{\ell-2\lambda}Q_{n+1}^{-\lambda}$ and corresponding highest weight $(\ell-2\lambda)\varpi_{n-1}$, whence $M^\lambda=L((\ell-2\lambda)\varpi_{n-1})$.
To show irreducibility, it remains to prove that the only primitive subspace of $M$ is $M^0$: to this purpose, we need to compute the action of $\mathfrak r_+^\ell$ on $M$, which is spanned by the differential operators $P_iP_jP_{n+1}$, $i,j=1,\dots,n$.
It suffices to show that there are no primitive subspaces except $M^0$ w.r.t.\ the action of $P_n^2P_{n+1}$: this is achieved by an easy computation.
In particular, we can show that the primitive monomial $Q_n^\ell$ generates the whole module $M$: a general element $Q^\mu$ of $H^0(\widetilde X ,\mathcal O(\ell D_0))$ is associated to a multiindex $\mu$ such that $\sum_{i=1}^n\mu_i-2\mu_{n+1}=\ell$.
The differential operator $D=Q^\mu P_n^\ell/\ell!$ belongs to $\mathcal D_{\mathcal O(\ell D_0)}(\widetilde X )$ by construction, and it is easy to verify that $Q^\mu$ is obtained from $Q_n^\ell$ by applying $D$.
\end{proof}
We consider now a non-positive integer $\ell$.
\begin{Thm}\label{t-0res-neg}
For any non-positive integer $\ell$, there is an isomorphism 
\[
H^0(\widetilde X ,\mathcal O(\ell D_0))\cong \begin{cases}
L_{\mathfrak{sp}_{2n}}(-\frac{1}2\varpi_n),& \text{if $\ell$ is even},\\
L_{\mathfrak{sp}_{2n}}(\varpi_{n-1}-\frac{3}2\varpi_n),& \text{if $\ell$ is odd},
\end{cases}
\]
of $\mathfrak g$-modules.
\end{Thm}
\begin{proof}
First of all, $H^0(\widetilde X ,\mathcal O(\ell D_0))$ is an $\mathfrak A_\ell$-module.
The weight decomposition of $M=H^0(\widetilde X ,\mathcal O(\ell D_0))$ corresponds, by the same arguments in the proof of Theorem~\ref{t-0res-pos}, to the decomposition of Theorem~\ref{t-cohres}:
\[
H^0(\widetilde X ,\mathcal O(\ell D_0))=\begin{cases}
\bigoplus_{\lambda\leq \frac{\ell}2}M^\lambda,& M^\lambda=H^0(\mathbb P^{n-1},\mathcal O(\ell-2\lambda)),\quad \text{$\ell$ even},\\
\bigoplus_{\lambda\leq \frac{\ell-1}2}M^\lambda,& M^\lambda=H^0(\mathbb P^{n-1},\mathcal O(\ell-2\lambda)),\quad \text{$\ell$ odd}.
\end{cases}
\] 
The factors in both direct sums are well-known to be irreducible, finite-dimensional highest-weight modules for $\mathfrak{sl}_n$, which can be identified with $L((\ell-2\lambda)\varpi_{n-1})$.
The infinite sequence of weights w.r.t.\ $\mathfrak z_\ell$ is entirely contained in the set of negative integers: hence, the requirements of Theorem~\ref{t-lift} are satisfied, and $M$ is a $\mathfrak g$-module.  

Using the results of Subsection~\ref{ss-1-3}, in particular (\ref{eq-chevgen}), the only primitive vectors of $H^0(\widetilde X ,\mathcal O(\ell D_0))$ are $(i)$ $Q_{n+1}^{-\frac{\ell}2}$, if $\ell$ is even, and $(ii)$ $Q_n Q_{n+1}^{-\frac{\ell-1}2}$, if $\ell$ is odd.
Moreover, both primitive vectors are weight vectors, with respective weights $-\frac{1}2\varpi_n$ and $\varpi_{n-1}-\frac{3}2\varpi_n$.
It remains to show that they generate the respective cohomology groups as $\mathfrak g$-modules: in fact, to any monomial $Q^\mu$ in $H^0(\widetilde X ,\mathcal O(\ell D_0))$, where $\mu$ is a multiindex satisfying $\sum_{i=1}^n\mu_i-2\mu_{n+1}=\ell$, we can associate the differential monomials $Q^\mu P_{n+1}^{-\frac{\ell}2}/(-\frac{\ell}2)!$, $\ell$ even, and $Q^\mu P_nP_{n+1}^{-\frac{\ell-1}2}/(-\frac{\ell-1}2)!$, $\ell$ odd.
It is easy to verify that both differential monomials are elements of $\mathcal D_{\mathcal O(\ell D_0)}(\widetilde X )$ in the respective cases.
Hence, the $0$-th cohomology is irreducible, being a highest weight vector with exactly one primitive vector in both cases $\ell$ even and odd, whence the claim.
\end{proof}
\begin{Rem}
Observe that, for $\ell$ even, there is an isomorphism of $\mathfrak g$-modules between the module $\mathcal O(\overline X)$ of regular functions on $\overline X$ of $\mathfrak g$ and the module $H^0(\widetilde X ,\mathcal O(\ell D_0))$, induced by mapping the highest weight vector $1$ of $\mathcal O(\overline X)$ to the highest weight vector $Q_{n+1}^{\frac{\ell}2}$ of $H^0(\widetilde X ,\mathcal O(\ell D_0))$; $\mathcal O(\overline X)$ was proved to be an irreducible highest weight $\mathfrak g$-module by other methods in~\cite{Levasseuretal}.
On the other hand, the highest weight of $H^0(\widetilde X ,\mathcal O(\ell D_0))$, for $\ell$ odd, belongs to the same Weyl orbit of the highest weight of $H^0(\widetilde X ,\mathcal O(\ell D_0))$, for $\ell$ even, as can be verified by an easy computation.
\end{Rem}
As for the $n-1$-th cohomology $H^{n-1}(\widetilde X ,\mathcal O(\ell D_0))$, we need only consider the case $\ell\leq -n$ (otherwise it is trivial).
\begin{Thm}\label{t-hres}
For any integer $\ell\leq -n$, there is a weight space decomposition
\[
H^{n-1}(\widetilde X ,\mathcal O(\ell D_0))\cong \bigoplus_{\lceil \frac{\ell+n}2\rceil\leq\lambda\leq 0} L_{\mathfrak{sl}_n}(-(\ell+n-2\lambda)\varpi_{1})
\]
in the category $\mathcal O_{\mathfrak A_\ell}$; the module $H^{n-1}(\widetilde X ,\mathcal O(\ell D_0))$ is irreducible and is of highest weight.
\end{Thm}
\begin{proof}
We first show that $H^{n-1}(\widetilde X ,\mathcal O(\ell D_0))$ belongs to the category $\mathcal O_{\mathfrak A_\ell}$: this is a consequence of Theorem~\ref{t-cohres}, namely
\[
M=H^{n-1}(\widetilde X ,\mathcal O(\ell D_0))=\bigoplus_{\lceil \frac{\ell+n}2\rceil\leq \lambda\leq 0} M^\lambda,\quad M^\lambda=H^{n-1}(\mathbb P^{n-1},\mathcal O(\ell-2\lambda)). 
\]
Every weight subspace $M^\lambda$ is an irreducible highest weight $\mathfrak{sl}_n$-module with highest weight vector $Q_1^{\ell+n-2\lambda-1}Q_2^{-1}\cdots Q_n^{-1}Q_{n+1}^{-\lambda}$ and corresponding highest weight $-(\ell+n-2\lambda)\varpi_1$, whence the identification $M^\lambda=L(-(\ell+n-2\lambda)\varpi_1)$; observe that the highest weight vector is a multiple of $e_{-\beta}^{-\lambda}(Q_1^{\ell+n-1}Q_2^{-1}\cdots Q_n^{-1})$, $\beta$ being the highest root of $\mathfrak g$.
Finally, $M$ is generated by the highest weight vector of $M^0$ as an $\mathfrak A_\ell$-module: in fact, any element of $H^{n-1}(\widetilde X ,\mathcal O(\ell D_0))$ has the form $Q^\mu$, where $\mu_i<0$, $i=1,\dots,n$, $\mu_{n+1}\geq 0$, and $\sum_{i=1}^n \mu_i-2\mu_{n+1}=\ell$. 
To such a monomial, we associate the differential operator 
\[
D_\mu=\frac{(-1)^{\sum_{i=1}^n \mu_i+n}}{\prod_{i=0}^{-\mu_1} (\ell+n-1-i)\prod_{j=2}^n (-(\mu_j+1))!}Q_1^{-\ell-n}Q_{n+1}^{\mu_{n+1}} P_1^{-\mu_1-1}\cdots P_n^{-(\mu_n+1)},
\]
which is readily checked to belong to $\mathcal D_{\mathcal O(\ell D_0)}(\widetilde X )$; it is also easy to verify that $Q^{\mu}=D_\mu (Q_1^{\ell+n-1}Q_2^{-1}\cdots Q_n^{-1})$.
\end{proof}

\subsubsection{The module structure on the cohomology of $Y$}\label{sss-4-2-2}
For any even integer $\ell$, we consider the cohomology of the weighted $n$-dimensional projective space $Y$ with values in the Serre bundle $\mathcal O(\ell)$, see Subsection~\ref{ss-3-2}.
As a byproduct of Corollary~\ref{cor-sympl}, the cohomology groups $H^\bullet(X,\mathcal O(\ell))$ inherit the structure of $\mathfrak A$-modules.
Recall the results of Subsection~\ref{ss-3-2} about the cohomology of $Y$ with values in $\mathcal O(\ell)$.
Recalling Theorem~\ref{t-weightres}, the Fourier transform $F_I$, for $I$ as in Subsection~\ref{ss-4-1} is used to compute the image of $\mathfrak A_{\ell+2}$ in $\mathcal D_{\mathcal O(\ell)}(Y)$ from $\mathcal D_{\mathcal O((\ell+2)D_0)}(\widetilde X )$: since the $I$-reflection affects only the generating vector $k_{n+1}'$, the only Chevalley generators of $\mathcal D_{\mathcal O((\ell+2)D_0)}(\widetilde X )$ affected by $F_I$ are  
\begin{align*}
F_I\!\left(\frac{1}2P_n^2 P_{n+1}\right)&=-\frac{1}2P_n^2 Q_{n+1},\quad &\quad F_I\!\left(-\frac{1}2 Q_n^2Q_{n+1}\right)&=\frac{1}2Q_n^2 P_{n+1},\\ 
F_I(-Q_{n+1}P_{n+1})&=Q_{n+1}P_{n+1}+1.
\end{align*}
Similarly to what was done in Subsubsection~\ref{sss-4-2-1}, we can identify the cohomology of $X$ with values in $\mathcal O(\ell)$, for any even integer $\ell$, with certain irreducible, highest weight $\mathfrak A_{\ell+2}$-modules in the category $\mathcal O_{\mathfrak A_{\ell+2}}$.
\begin{Thm}\label{t-0weight}
For any even, non-negative integer $\ell$, there is a weight space decomposition
\[
H^0(Y,\mathcal O(\ell))\cong \bigoplus_{1\leq \lambda \leq \frac{\ell+2}2} L_{\mathfrak{sl}_n}((\ell-2\lambda+2)\varpi_{n-1})
\]
in the category $\mathcal O_{\mathfrak A_{\ell+2}}$; the module $H^0(Y,\mathcal O(\ell))$ is irreducible and of highest weight.
\end{Thm}
\begin{proof}
The global sections of $\mathcal O(\ell)$ are spanned by monomials $Q^\mu$, with multiindices $\mu$ satisfying $\sum_{i=1}^n\mu_i+2\mu_{n+1}=\ell$: hence, $0\leq \mu_{n+1}\leq \frac{\ell}2$.
The weight subspaces $M^\lambda$ of $M=H^0(Y,\mathcal O(\ell))$ correspond to weights $1\leq \lambda\leq \frac{\ell+2}2$: in fact, $M$ splits as $M=\bigoplus_{1\leq m\leq \frac{\ell+2}2} M^\lambda$, with $M^\lambda$ the space of complex polynomials in $n$ variables, homogeneous of degree $\ell-2\lambda+2$, the isomorphism being concretely 
\[
M^\lambda \ni Q_1^{\mu_1}\cdots Q_n^{\mu_n}\mapsto Q_1^{\mu_1}\cdots Q_n^{\mu_n}Q_{n+1}^{\lambda-1}\in M.
\]   
Every weight subspace $M^\lambda$ is an irreducible, finite-dimensional highest weight module for $\mathfrak{sl}_n$: the highest weight vector is $Q_n^{\ell-2\lambda +2}Q_{n+1}^{\lambda-1}$ with highest weight $(\ell-2\lambda+2)\varpi_{n-1}$, whence the identification $M^\lambda=L((\ell-2\lambda+2)\varpi_{n-1})$.
Observe that the $\mathfrak{sl}_n$-highest weight vector of $M^\lambda$ is a multiple of $e_{-\alpha_n}^{\ell-2\lambda+2}(Q_{n+1}^{\frac{\ell}2})$, $0\leq 0\leq \frac{\ell}2$, with $Q_{n+1}^{\frac{\ell}2}$ the highest weight vector of $M^{\frac{\ell+2}2}$.
Therefore, $H^0(Y,\mathcal O(\ell))$ belongs to the category $\mathcal O_{\mathfrak A_{\ell+2}}$; its highest weight w.r.t.\ $\mathfrak z_{\ell+2}$ belongs to $\mathbb N$, the sequence of corresponding weights has finitely many terms, hence Theorem~\ref{t-lift} does not apply.

There is only one primitive one-dimensional subspace, spanned by the vector $Q_{n+1}^{\frac{\ell}2}$ in $H^0(Y,\mathcal O(\ell))$; hence, the module is irreducible.
The generator $Q_{n+1}^{\frac{\ell}2}$ is a highest weight vector, which additionally generates $H^0(Y,\mathcal O(\ell))$ as an $\mathcal O_{\mathfrak A_{\ell+2}}$-module by computations similar to those in the proof of Theorem~\ref{t-0res-neg}.
\end{proof}
As for the $n$-th cohomology of $Y$ with values in $\mathcal O(\ell)$, for an even integer $\ell$, we have
\begin{Thm}\label{t-nweight}
If $\ell$ is any even integer, such that $\ell\leq -n-2$, there is a weight space decomposition
\[
H^n(Y,\mathcal O(\ell))\cong \bigoplus_{\lceil \frac{\ell+n}n\rceil \leq\lambda\leq 0} L_{\mathfrak{sl}_n}(-(\ell+n-2\lambda+2)\varpi_1)
\]
of $\mathfrak A_{\ell+2}$-modules; furthermore, $H^n(Y,\mathcal O(\ell))$ is irreducible and is generated by a highest weight vector.
\end{Thm}
\begin{proof}
We recall that $M=H^n(Y,\mathcal O(\ell))$ is spanned by negative monomials $Q^\mu$, with $\mu_i<0$, $i=1,\dots,n+1$, and $\sum_{i=1}^n\mu_i+2\mu_{n+1}=\ell$, thus, we have the splitting  
\[
M=\bigoplus_{\lceil\frac{\ell+n}2\rceil\leq m\leq -1} S^n_m,
\]
where $S_m^n$ denotes the subspace of negative monomials in $n$ variables of the form $Q_1^{\mu_1}\cdots Q_n^{\mu_n}$ with $\mu_i<0$ and $\sum_{i=1}^n\mu_i=\ell-2m$, the isomorphism being given by
\[
S_m^n\ni Q^{\mu_1}\cdots Q_n^{\mu_n}\mapsto Q_1^{\mu_1}\cdots Q_n^{\mu_n}Q_{n+1}^m\in M.
\]
(Observe that the negative monomials must be considered modulo coboundaries.)
The subspace $S_m^n$ is an irreducible highest weight module for $\mathfrak{sl}_n$: the highest weight vector is $Q_1^{\ell+n-2m-1}Q_2^{-1}\cdots Q_n^{-1}$ with highest weight $-(\ell+n-2m)\varpi_1$, hence, recalling the $\mathfrak z_{\ell+2}$-action on $M$, $M^\lambda=L(-(\ell+n-2\lambda+2)\varpi_1)$, $\lambda=m+1$. 
The module $M$ belongs to the category $\mathcal O_{\mathfrak A_{\ell+2}}$; by the previous arguments it is finite-dimensional; the only primitive subspace of $M$ corresponds, by direct computations, to $M^0$, whence $M$ is irreducible.
Finally, the highest weight vector of $M^0$ generates $M$ as an $\mathfrak A_{\ell+2}$-module by slight modifications of the arguments in the proof of Theorem~\ref{t-hres}.  
\end{proof}

\end{document}